\documentclass{amsart}
\usepackage{latexsym}
\usepackage{amssymb}
\usepackage{enumerate}
\usepackage[dvips]{graphicx}
\usepackage[latin1]{inputenc}
\title[Hyperbolic and conformal measures]{Measures with positive Lyapunov exponent and conformal measures in rational dynamics }
\author{Neil Dobbs}
\address{KTH, Stockholm} 
\email{neil.dobbs@gmail.com}

\newcommand\beginpf{\begin{proof}}
\newcommand\eprf{\end{proof}}

\newcommand\Crit{\mathrm{Crit}}
\newcommand\arr{\mathbb{R}}

\newcommand\dist{\mathrm{dist}}
\newcommand\Fix{\mathrm{Fix}}
\newcommand\cal{\mathcal}

\newcommand\scrM{{\cal M}}
\newcommand\scrP{{\cal P}}

\newcommand\spher{\sigma}
\newcommand\omu{{\overline{\mu}}}
\newcommand\J{{\cal{J}}}
\newcommand\HD{\mathrm{HD}}

\newcommand\cbar{{\overline{\mathbb{C}}}}
\newcommand\N{\mathbb{N}}
\newcommand\remark{\noindent \emph{Remark: }}

\newtheorem{thm}{Theorem}
\newtheorem{dfn}[thm]{Definition}
\newtheorem{lem}[thm]{Lemma}
\newtheorem{prop}[thm]{Proposition}
\newtheorem{cor}[thm]{Corollary}
\newtheorem{Fact}[thm]{Fact}

\begin{document}
\date{\today}
\maketitle
\begin{abstract}
Ergodic properties of rational maps are studied, generalising the work of F.\ Ledrappier. 
A new construction allows for simpler proofs of stronger results.
Very general conformal measures are considered. Equivalent conditions are given for 
 an ergodic invariant probability measure with positive Lyapunov exponent to be absolutely continuous with respect to a general conformal measure.  If they hold, we can construct an induced expanding Markov map with integrable return time which generates the invariant measure. 
\end{abstract}

\section{Introduction}
\subsection{Overview}
For dynamical systems defined by rational maps of the Riemann sphere, we study ergodic invariant measures with positive Lyapunov exponents, generalising and improving upon the important work of Ledrappier (\cite{Ledrappier:AbsCnsComplex}). We construct a set with good \emph{Markov} properties in the natural extension. In Theorem \ref{thm:dbbdense} we give equivalent properties for measures to be absolutely continuous with respect to a general conformal measure. 

In the presence of chaotic dynamics, ergodic invariant probability measures give a good statistical description of the dynamics of points seen by the measure. In general there are many such measures and one is interested in measures which see a large collection of points. Given a reference measure such as Lebesgue measure (if the space is some subset of Euclidean space), one could look for a measure which is absolutely continuous with respect to the reference measure. 

In complex dynamics, the Julia set of a rational map $f : \cbar \to \cbar$ is often a fractal set of zero Lebesgue measure. When this happens the notion of a $t$-conformal measure generalises that of Lebesgue measure (which is a 2-conformal measure). 
A $t$-conformal measure always exists (\cite{Sullivan:Conformal}) for some $t>0$. Its support is the whole Julia set and it is the useful reference measure. In situations with sufficient expansion 
 (\cite{PRL:StatTCE, GraczykSmirnov} amongst others) there exists an ergodic invariant probability measure which is absolutely continuous with respect to the (in this case unique) $t$-conformal measure. 

When an ergodic absolutely continuous invariant probability measure (\emph{acip}) with respect to a $t$-conformal measure can be shown to exist, the construction can be used to show  the measures are equivalent and that the acips have density bounded from below away from zero and positive Lyapunov exponent. We show that whenever such an acip with positive exponent exists, there exists an expanding induced Markov map with integrable return time which generates the measure (see Definition \ref{dfn:markov}) and the density is bounded away from zero. 

More general definitions of conformal measures exist (see Definition \ref{dfn:phitconforme}), and are related to equilibrium measures for certain potential functions. Our Theorem \ref{thm:dbbdense} below states that, given a conformal measure and an ergodic invariant probability measure $\mu$ with positive Lyapunov exponent, there is equivalence between: $\mu$ being absolutely continuous; $\mu$ satisfying a variant of Pesin's formula; the existence of an expanding induced map which generates $\mu$. Moreover, when these conditions are satisfied the density of $\mu$ with respect to the conformal measure is  bounded away from zero.  

An open set $U$ is called \emph{nice}, or \emph{regularly returning}, if $f^n(\partial U) \cap U = \emptyset$ for all positive $n$. This condition implies that if $A$ and $B$ are connected components of $f^{-n_a}(U)$ and $f^{-n_b}(U)$ respectively, for some $n_a,n_b \geq 0$, then $A$ and $B$ are either \emph{nested or disjoint}, that is, $A \subset B$ or $B \subset A$ or $A \cap B = \emptyset$. 
For polynomial maps the well-known Yoccoz \emph{puzzle pieces} are nice sets. Nice sets have also been constructed by Rivera-Letelier (\cite{Rivera:Connecting}) for backward-contracting rational maps, and used again by Przytycki and Rivera-Letelier in \cite{PRL:StatTCE}. However, for more general rational maps it seems unlikely that nice sets exist. 

Therefore, in order to get the good topological properties of an induced Markov map, we need to select our inverse branches (looking backwards, or the iterates looking forward) with care. 
The existence of the induced Markov map is new and is based on an
 interesting new concept introduced in this paper:  the \emph{regularly returning cylinder} (Definition \ref{dfn:regretcyl}) in the natural extension.
The regularly returning cylinder in the natural extension will have this \emph{nested or disjoint} property for projected preimages $\emph{without}$ being a nice set itself. 

This work was inspired by the results of Ledrappier stated in \cite{Ledrappier:AbsCnsComplex}. We give a very different proof of the existence of the unstable manifold in the natural extension. The regularly returning cylinders we construct allow us to show (for the measure considered) the existence of a \emph{finite} generating partition with good topological properties. This is completely new as far as we are aware. Since the partition constructed has one more element than the degree $d$ of the rational map, it is not hard to deduce that the entropy of the measure is bounded by $\log d$.  This gives an alternate way of proving that the topological entropy is $\log d$. 

A related construction can be found in the historically important article of Mañé, \cite{Mane:Bernoulli}. He proved the existence of a good collection of inverse branches of some iterate $f^m$ of $f$ for the measure of maximal entropy, where $m > 1$ depends on $f$. Denker and Urbanski (\cite{denker1991ete}) extended the construction to other measures with positive entropy (and thus positive Lyapunov exponent), but $m$ then also depends on the measure. The technique of Mañé relies heavily on the one-dimensional holomorphic setting. Our construction, on the other hand, is based on Pesin theory and more easily generalisable. It also gives a finite partition with different and interesting properties.

Passing on to the study of absolutely continuous measures, again our proofs are relatively simple thanks to the regularly returning cylinders. Related 
results for equilibrium measures with respect to  Hölder continuous potentials can be found in \cite{przytycki1990pfr, denker1991ete}. The singularity of $\log|Df|$, due to the presence of 
critical points, complicates matters when one wants to consider more interesting potentials. For the 
potential $-t_0 \log|Df|$, where $t_0$ is the hyperbolic dimension of the Julia set or, equivalently, the minimal exponent for which a $t$-conformal measure exists, \cite{Ledrappier:AbsCnsComplex} contains important results. We extend these to show that existence of an acip implies the existence of an  induced Markov map with integrable return time, and give an explicit proof that the density is bounded away from zero, a fact used in \cite{PRL:StatTCE}. 

Our definitions of regularly returning cylinders and of Markov maps and our proof of the existence of the unstable manifold are not stated in the language of rational dynamics in order to show that these ideas are not particular to holomorphic dynamics and will generalise to a higher-dimensional non-conformal setting. We nevertheless include extendibility in the definitions for those, including the author, who like their ``Koebe space''.  

\subsection{Definitions and main result}
As usual when dealing with rational maps, we make use of the spherical metric and derivative. We denote the associated spherical distance by $\spher(\cdot, \cdot)$ and the spherical derivative at a point $x$ by $|Df(x)|$. We use $|\cdot|$ to denote the diameter of a set. All subsets of $\cbar$ considered in this paper are assumed to be Borelian.

\begin{dfn} A measure $m$ is \emph{$t$-conformal} for rational  map $f: \cbar \to \cbar$ and $t \in \arr$ if $m(\J) = m(\cbar) = 1$ and 
$$
m(f(A)) = \int_A |Df|^t dm
$$
 for any set $A$ on which $f$ is injective. 
\end{dfn}
We shall work with a more general definition. 
Let $\phi : \cbar \to \cbar$ be a H\"older continuous function and let $t \in \arr$. 
Consider the relation
\begin{equation}\label{eqn:conformality}
m(f(A)) = \int_A e^\phi |Df|^t dm
\end{equation}
\begin{dfn} \label{dfn:phitconforme}
If $t \geq 0$, a measure $m$ is \emph{$(\phi,t)$-conformal} for a rational map $f$ if $m(\J) = m(\cbar) = 1$, if (\ref{eqn:conformality}) holds for
 every   set $A$ on which $f$ is injective. 
If $t < 0$, a measure $m$ is \emph{$(\phi,t)$-conformal} if $m(\J) = m(\cbar) = 1$, if (\ref{eqn:conformality}) holds for
 every  set $A \subset \cbar \setminus \Crit$ on which $f$ is injective  and if $m(Crit) = 0$ (see also the related definition of \emph{backward conformal measures} in \cite{PrzRLSm}, similar when $\phi$ is constant). 
\end{dfn}
Spherical Lebesgue measure is a $(0,2)$-conformal measure. For rational maps of degree $d\geq 2$ the measure of maximal entropy has constant Jacobian equal to $d$ so it is a $(\log d,0)$-conformal measure. 
\begin{dfn} We denote by $\scrM(f)$ the collection of ergodic, $f$-invariant probability measures.
\end{dfn}
\begin{dfn} The \emph{Lyapunov exponent} of a measure $\mu \in \scrM(f)$ is defined as $\chi_\mu := \int \log |Df| d\mu$. The entropy of $\mu$ we denote by $h_\mu$.
\end{dfn}
Of course $-\infty \leq \chi_\mu < \infty$.
If $\mu$ is not supported on a periodic attractor, then $\chi_\mu \geq 0$ by \cite{Przytycki:LyapNonneg}.
\begin{dfn} We denote by $\scrM^+(f)$ the collection of $\mu \in \scrM(f)$ for which $\chi_\mu > 0$.
\end{dfn}
\begin{dfn} We call a $(\phi,t)$-conformal measure $m$  \emph{exceptional} if $\mathrm{Supp}(m) \ne \J$, and call it \emph{non-exceptional} if it admits no restriction which, when normalised, is an exceptional conformal measure.
\end{dfn}
\remark If  $t \geq 0$ then $m$ is non-exceptional by the \emph{eventually onto} property of rational maps. Moreover, unless we are in a special case where there is a finite forward-invariant set $E$ contained in the  post-critcal set $PC := \bigcup_{i\geq 0} f^i(f(\Crit))$ such that 
$$
f^{-1}(E) \subset E\cup \Crit,$$
then even if $t < 0$ the measure $m$ is necessarily non-exceptional. If $m$ is exceptional, $m$ is atomic and is supported on such a set $E$. This exceptional setting is very rare; interesting examples are Lattès maps.

\remark In the absence of parabolic periodic points in $PC$, exceptional and non-exceptional $(\phi,t)$-conformal measures cannot co-exist for a given $(\phi, t)$. 
\begin{dfn} \label{dfn:markov} Let $f$ be a rational map with a $(\phi,t)$-conformal measure $m$. Let $U$ and $V$ be open topological balls such that $\overline{U} \subset V$ and 
 $m(U) > 0$ and let $U_i \subset V_i \subset U$ satisfy:
\begin{itemize}
\item $U_i \cap U_j \ne \emptyset \iff i = j$;
\item $m(U \setminus \bigcup_i U_i) = 0$;
\item for each $i$ there exists $n_i$ such that $f^{n_i} : V_i \to V$ is a diffeomorphism with $|Df^{n_i}| > 2$ and $f^{n_i}(U_i) = U$;
\item
    there exist $C, \delta > 0$ such that, for each $i$ and all $j\leq n_i $, 
    $$
    |Df^j(x)| > C e^{\delta j}
    $$
    for all $x$ in  $f^{n_i - j}(U_i)$. 
\end{itemize}
Define $\psi: \bigcup_i U_i \to U$ by $\psi_{|U_i} := f^{n_i}_{|U_i}$. Then $\psi$ is called an \emph{expanding induced Markov map for $(f,m)$}. If additionally
 $\sum_i n_i m(U_i) < \infty$ then we say $\psi$ has \emph{integrable return time}. Let $\nu$ be the unique (c.f. Lemma \ref{lem:markovac}) absolutely continuous $\psi$-invariant probability measure. It has bounded density. The following measure, when normalised, is a well-defined, ergodic, absolutely continuous $f$-invariant probability measure which we say is \emph{generated by $\psi$}: 
$$
\overline{\nu} := \sum_i \sum_{j=0}^{n_i -1} f_*^j \nu.
$$
\end{dfn}

\begin{thm} \label{thm:dbbdense} Let $f$ be a rational map of the Riemann sphere and $m$ a non-exceptional $(\phi,t)$-conformal measure for $f$. Let $\mu \in \scrM^+(f)$. The following conditions are  equivalent:
\begin{enumerate}[$1^\circ$.]
\item \label{itm:ac} $\mu <\!< m$;
\item \label{itm:dbb} there exist $\varepsilon >0$ and a measurable density function $\rho \geq \varepsilon > 0$ such that $\mu = \rho m$ on an open set $U$ of positive measure, and if $t \geq 0$ one can take $U = \cbar$;
\item \label{itm:pesin} $h_\mu = t\chi_\mu + \int \phi d\mu$;
\item \label{itm:hd} $\HD(\mu) = t + \frac{\int \phi d\mu}{\chi_\mu}$;
\item\label{itm:mark}  
 there exists an expanding induced Markov map for $(f,m)$ with integrable return time which generates $\mu$.
\end{enumerate}
Should the equivalent conditions hold then $m$ is  ergodic and unique and $\mu$ is the unique measure in $\scrM^+(f)$ satisfying any of the equivalent conditions. 
\end{thm}
\beginpf Clearly \ref{itm:mark}$^\circ$ $\Rightarrow$ \ref{itm:ac}$^\circ$. 
In Proposition \ref{prop:complexHDmu} we show \ref{itm:ac}$^\circ \Rightarrow$ \ref{itm:hd}$^\circ$.
By the Dynamical Volume Lemma (Proposition \ref{prop:ratdvl}) \ref{itm:hd}$^\circ \iff$ \ref{itm:pesin}$^\circ$. 
In Proposition \ref{prop:acmarkov} we prove \ref{itm:pesin}$^\circ$ $\Rightarrow$ \ref{itm:dbb}$^\circ$ and  \ref{itm:pesin}$^\circ$ $\Rightarrow$ \ref{itm:mark}$^\circ$. 
\eprf

\remark If $m$ is exceptional then one can show easily that \ref{itm:ac}$^\circ$ implies \ref{itm:pesin}$^\circ$, \ref{itm:hd}$^\circ$ and \ref{itm:mark}$^\circ$, 
and that \ref{itm:mark}$^\circ \Rightarrow$ \ref{itm:ac}$^\circ$. We no longer necessarily have that \ref{itm:pesin}$^\circ$ implies absolute continuity: 
consider the quadratic Chebyshev  map $x \mapsto 4x(1-x)$. For $t = -1$, there are two equilibrium measures. One is the measure for maximal entropy, the other is the unit mass on $\{0\}$ which is absolutely continuous with respect to the unique $(\log 4, -1)$-conformal measure. 

\subsection{Equilibrium measures and the pressure function}
In carrying out this work, we bore in mind the idea that a measure $\mu$ satisfying the equivalent conditions in the preceding theorem should be the equilibrium measure for some potential function related to $\phi$ and to $t$. 

One can consider the (variational) pressure
$$
P(\phi,t) := P_{\mathrm{var}}(\phi,t) := \sup_{\mu \in \scrM(f)}(h_\mu + \int (-\phi -t \log|Df|) \, d\mu),
$$
(for the potential $-\phi -t\log|Df|$). A measure which realises the supremum is called an equilibrium measure for that potential. A measure $\mu$ verifying the equivalent conditions in  Theorem \ref{thm:dbbdense} is an equilibrium measure provided  $P(\phi, t) = 0$. 

Given some function $\phi'$, one can set $\phi := \phi' + P(\phi',t)$. Then $P(\phi, t) = 0$. 
 One can expect that, when $P(\phi,t) = 0$, there will exist a $(\phi, t)$-conformal measure. When $t = 0$ this can be found in \cite{denker1991ete, przytycki1990pfr};  existence of equilibrium states   is also shown. 

Now let us consider the 
 \emph{natural} potential $-t \log |Df|$ and the corresponding pressure function 
$$
P(t) := P(0,t) = \sup_{\mu \in \scrM(f)} h_\mu - \int t \log |Df| \, d\mu .
$$
Przytycki, Rivera-Letelier and Smirnov, generalising the work of Denker and Urba\'nski and Przytycki (\cite{DenkerUrbanski:RatConf, Przytycki:LyapNonneg}), have shown the existence of $(P(t),t)$-conformal measures for rational maps for all $t \in \arr$ (proof of proposition 1.2 and section A.2 in \cite{PrzRLSm}). 
Moreover, they show that for $t \geq 0$, $P(t)$ is the infimum of $p \in \arr$ such that there exists a $(p,t)$-conformal measure. For $t < 0$, if $m$ is a non-exceptional $(p,t)$-conformal measure, they show an even stronger statement, $p = P(t)$. 
\begin{cor}
Let $f$ be a rational map of the Riemann sphere, $t \geq 0$, and let $\mu \in \scrM^+(f)$. If $p \in \arr$ and $m$ is a $(p,t)$-conformal measure such that $\mu <\!< m$, then $p = P(t)$ and $\mu$ is an equilibrium measure.  
\end{cor}
\beginpf We know already $p \geq P(t)$. By Theorem \ref{thm:dbbdense}, $p = h_\mu - t\chi_\mu$, thus $p \leq P(t)$. 
\eprf
\begin{cor}
Let $f$ be a rational map of the Riemann sphere, $t \geq 0$, let $m$ be a $(P(t),t)$-conformal measure and let $\mu \in \scrM^+(f)$. Suppose $h_\mu - t\chi_\mu = P(t)$. Then $\mu$ is equivalent to $m$. In particular, $\mu$ is the unique equilibrium state with positive Lyapunov exponent for the potential $-t \log |Df|$ and $m$ is ergodic. 
\end{cor}
\remark Przytycki (\cite{Prz06:Bull}) conjectured that, for $t$ equal to the first zero of the pressure function, there was a unique equilibrium state with positive exponent. For that particular value of $t$, the answer can be shown using the work of Ledrappier (\cite{Ledrappier:AbsCnsComplex}). Rephrasing the previous corollary, we answer a stronger question:
\begin{cor}
For each $t \geq 0$, 
there is at most one equilibrium state with positive Lyapunov exponent for the potential $-t \log |Df|$.
\end{cor}
\remark Provided $P(t) \ne 0$, any equilibrium state for the potential $-t \log |Df|$ has positive Lyapunov exponent. 

The case $t < 0$ was studied by Makarov and Smirnov (\cite{MakarovSmirnov:Nonrec}) who used operator theory to show existence of equilibrium states and conformal measures and to study uniqueness of equilibrium states and analyticity of the pressure function. One can approach this differently. As recalled in \cite{PrzRLSm}, the Patterson-Sullivan method for constructing conformal measures easily produces a $(P(t),t)$-conformal measure (one does need to take care with the definitions of pressure). 
Let 
$$
t_0:= \inf \{t : \mbox{ there exists a non-exceptional $(P(t),t)$-conformal measure} \}.$$
\begin{prop} For each $t \in (t_0,0]$ there is a unique equilibrium measure and a unique $(P(t),t)$-conformal measure.
\end{prop}
\beginpf
Existence of equilibrium measures follows from upper semi-continuity of $h_\mu$ and $\chi_\mu$. Let $m$ be a $(P(t),t)$-conformal measure, necessarily non-exceptional since exceptional and non-exceptional measures cannot co-exist. By Theorem~\ref{thm:dbbdense}, $m$ and $\mu$ are equivalent, thus $m$ is ergodic. Any other equilibrium state must also be equivalent to $m$ and thus to $\mu$ and thus coincides with $\mu$.
\eprf

For $t > 0$,
in a work-in-progress, Przytycki and Rivera-Letelier show the existence of equilibrium states for all $t$ in a neighbourhood of $[0,P^{-1}(0)]$ for a large class of maps (the \emph{Topological Collet-Eckmann} class). We expect existence of equilibrium states to hold for all rational maps for all $t$ strictly less than the first zero of the pressure function.

\subsection{Acknowledgements} This work grew out of the author's doctoral thesis, carried out at Université Paris-Sud (Orsay) under the supervision of Jacek Graczyk. I would like to thank Feliks Przytycki, Juan Rivera-Letelier and Mike Todd for helpful comments and discussions and the referee for helpful remarks. The exposition was substantially improved following a mini-course on the subject given by the author at Universidad Cat\'olica del Norte, Antofagasta; the attendees asked many helpful questions.  The author has been supported by the G\"oran Gustafssons Styftelse, the Knut och Alice Wallenbergs Styftelse and the EU Research and Training Network ``Conformal Structures and Dynamics'' (CODY). The author was also partially supported by the Research Network on Low Dimensional Dynamics, PBCT ACT 17, CONICYT, Chile.

\section{Unstable manifold}\label{sec:complexledrap}
\begin{dfn}{\emph{The Natural Extension \cite{Rohlin:Exact, Ledrappier:AbsCnsComplex}:}}
Let $Y \subset \cbar^\N$ be the set of sequences $y = (y_0,y_1,\ldots)$ such that $f(y_{i+1}) = y_i$  
and define $F: Y \to Y$ by
$$
F (y_0, y_1, \ldots) = (f(y_0),y_0,y_1,\ldots).
$$ 
Then $F : Y \to Y$ is invertible. We shall denote by $\Pi$ the projection onto the first element, i.e. $\Pi y = y_0$.

If $f$ has an invariant measure $\mu$ then there is a unique measure $\overline{\mu}$ on $Y$, invariant by $F$, such that $\Pi_* \overline{\mu} = \mu$. If $\mu$ is ergodic then so is $\overline{\mu}$, and importantly $\overline{\mu}$ is also an ergodic invariant measure for the inverse map $F^{-1}$.
We call $(Y, F, \overline{\mu})$ the natural extension of $(\cbar, f, \mu)$.
\end{dfn}


In this section we shall provide an alternative proof of the well-known result of Ledrappier (\cite{Ledrappier:AbsCnsComplex}) on existence of the unstable manifold in the natural extension. Our proof is very elementary and does not rely on non-flatness of the critical points. Recall that the spherical distance we denote by $\spher(\cdot,\cdot)$ and the spherical derivative of $f$ by $|Df|$.

\begin{prop}\label{prop:complexNeilUnstableManif}
Let $f$ be a rational map of $\cbar$. Suppose $f$ has an invariant, ergodic probability measure $\mu$ with positive Lyapunov exponent $\chi = \int \log|Df| \, d\mu > 0.$ Let $(Y, F, \overline{\mu})$ be the natural extension for $(f, \mu)$. 

Then there exists a measurable function $\alpha$ on $Y$, positive almost everywhere, such that, if $y \in Y$ and $\alpha(y) >0$:
\begin{itemize}
    \item
there exists a set $V_y \subset Y$ with  
$y \in V_y $ and $\Pi V_y = B(\Pi y, \alpha(y))$;
\item 
    for all $n \geq 0$, 
$f^n : \Pi F^{-n}V_y \to \Pi V_y$ is a diffeomorphism;
\item 
    for all $y' \in V_y$
$$
 \sum_{i=1}^\infty \left| \log|Df(\Pi F^{-i}y')| - \log |Df(\Pi F^{-i}y)|\right|  
< \log 2.
$$
\end{itemize}
For each $\eta >0$ there exists a measurable function $\rho$ on $Y$, positive almost everywhere, such that if $\rho(y) >0$ and if $n \geq 0$ and  $y' \in V_y$,  then
$$
\rho(y)^{-1}e^{n(\chi - \eta)} \leq |Df^n(\Pi F^{-n}y')| \leq \rho(y)e^{n(\chi + \eta)}.
$$
In particular, $|\Pi F^{-n}V_y| \leq \rho(y)e^{-n(\chi-\eta)} |\cbar|$.
\end{prop}

The strategy is as follows. First we show that we have a good distortion bound on $B(x,\min(c_0,|Df(x)|^3))$ for some constant $c_0$. Then we show that the derivative $|Df(\Pi F^{-n}y)|$ along backwards orbits is bounded from below by a sub-exponential sequence almost everywhere.  This allows us to define a sequence of balls on which one has (slow-) exponentially good distortion. Positive Lyapunov exponent will then imply that the pullbacks of some small ball will always land inside the balls with exponentially good distortion bounds, so the total distortion will be summable.

\begin{lem} \label{lem:ctwoproperty}
Let $f : \cbar \to \cbar$ be a rational map. Then there exists a $c_0 > 0$ with the following property. For all $c$ such that $0 < c < c_0$, if $x,x' \in \cbar$, $\spher(x,x') < c^3$ and $|Df(x)| > c$, then 
$$
\left| \log \frac{|Df(x)|}{|Df(x')|}\right| < 
c.$$
\end{lem}
\beginpf As a holomorphic map of the Riemann sphere, there exists a $C > 1$ such that for all $x,x' \in \cbar$, $|Df(x)| - |Df(x')| \leq C\sigma(x,x')$.
Dividing by $|Df(x)|$ and rearranging, it follows that
\begin{equation}\label{eqn:siggy}
1 -  \frac{C}{|Df(x)|} \sigma(x,x') \leq \frac{|Df(x')|}{|Df(x)|}  \leq 1+ \frac{C}{|Df(x)|} \sigma(x,x').
\end{equation}
Choose a $c_0 \in (0,1)$ small enough that $c_0C < 1$ and for all $r \in (0,1)$ we have $\log(1-rc_0C) > -r$. Suppose now that $|Df(x)| > c$ and $\sigma(x,x') < c^3$. Then 
$$
\frac{C}{|Df(x)|} \sigma(x,x') < c^2C 
< (c_0 C) c.$$
Taking logs of (\ref{eqn:siggy}),
$$
-c \leq \left| \log |Df(x)| - \log |Df(x')| \right| \leq c
$$
as required.
\eprf

We will need to swallow up some constants.  Fix $\delta > 0$ such that $\chi - \delta > 6\delta$ and $N > 0$ large enough that, for all $n \geq N$,  the following inequalities hold:
\begin{eqnarray}
\nonumber \label{eqn:sumndeltab} \frac{1}{2} \log 2 +\sum_{i\geq N}  e^{-i\delta} &<& \log 2;\\
\label{eqn:ballsb} 
2 e^{-n(\chi-\delta)} &<& 2^{-1}e^{-(n+1)4\delta };\\
\nonumber e^{-n\delta} &<& c_0.
\end{eqnarray}

Now, for each $y \in Y$, set $n(y)$ as the minimal $n \geq N$ such that 
    for all $n \geq n(y)$,
    \begin{equation}
        \label{eqn:n1}
|Df(\Pi F^{-n}y)| \leq 2 e^{-n\delta}
    \end{equation}
and 
    \begin{equation}
        \label{eqn:n2}
 e^{n(\chi - \delta)} \leq |Df^n(\Pi F^{-n}y)| \leq e^{n(\chi + \delta)}
    \end{equation}
 provided such an $n$ exists; otherwise set $n(y) := 0$. 
 Then $y \mapsto n(y)$ is a measurable function. 
\begin{lem} \label{lem:sea} For $\omu$ almost every $y$, $n(y) > 0$.
\end{lem}
\beginpf Inequality (\ref{eqn:n1}) holds almost everywhere for all large $n$ because the limit of $(1/n)\log |Df^n(\Pi F^{-n}y)|$ exists for almost every $y$; inequality (\ref{eqn:n2}) because the limit equals $\chi$ almost everywhere. 
\eprf

\begin{lem} \label{lem:bndistn}
    Suppose $n(y) > 0$. 
 Let $B_n := B(\Pi F^{-n}y, 2^{-1} e^{-n3\delta})$. For all $n \geq n(y)$, for all $x,x' \in B_n$,
$$
\left| \log|Df(x)| - \log|Df(x')|\right| < e^{-n\delta}.$$
\end{lem}
\beginpf
Follows from Lemmas \ref{lem:ctwoproperty}, \ref{lem:sea}.
\eprf
\begin{lem} 
    Suppose $n(y) > 0$. 
    For $n \geq n(y)$, $f(B_{n+1}) \supset B(\Pi F^{-n}y, 2 e^{n(\chi - \delta)})$.
\end{lem}
\beginpf By the preceding lemmas, $|Df(x)| > e^{-(n+1)\delta}$ on $B_{n+1}$, so $f(B_{n+1}) \supset B(\Pi F^{-n}y, 2^{-1}e^{-(n+1)3\delta}e^{-(n+1)\delta})$. Then use (\ref{eqn:ballsb}).
\eprf

\begin{lem} \label{lem:ballinduct}
    Suppose $n(y) > 0$. 
 Suppose $n \geq n(y)$ and $V$ is an open ball containing $\Pi y$ with $|V| \leq 1$ and suppose $V_n \subset B_n$ is such that $V_n \ni \Pi F^{-n}y$ and $f^n: V_n \to V$ is a diffeomorphism with distortion bounded by some $r$ with $0 < r < \log 2$, i.e.,
$$
\left| \log |Df^n(x)| - \log |Df^n(x')|\right| < r < \log 2$$
for all $x,x' \in V_n$. 
Then there exists $V_{n+1} \ni \Pi F^{-(n+1)}y$ such that the map $f^{n+1} : V_{n+1} \to V$ is a diffeomorphism with distortion bounded by 
$$
 r + e^{-(n+1)\delta}.$$ 
\end{lem}
\beginpf 
We have that $|Df^n(\Pi F^{-n}y)| > e^{n(\chi - \delta)}$ so $|V_n| < e^{-n(\chi - \delta)}$ and $V_n \subset f(B_{n+1})$. The result follows.
\eprf
\bigskip

\noindent
\begin{proof}[Proof of Proposition \ref{prop:complexNeilUnstableManif}]
If $n(y) >0$, 
let $V$ be the maximal  ball of radius bounded by $1$ and centred on $\Pi y$ such that there exists a set $V_{n(y)} \ni \Pi F^{-n(y)}(y)$ such that $f^{n(y)} : V_{n(y)} \to V$ is a diffeomorphism and such that for all $x,x' \in V_{n(y)}$, 
$$
\sum_{i=0}^{n(y)-1} \left| \log |Df(f^i(x))| - \log |Df(f^i(x'))|\right| \leq (1/2) \log 2.$$

For $0 \leq n < n(y)$ define $V_n := f^{n(y)-n}(V_{n(y)})$. For $n > n(y)$ define $V_n$ inductively using Lemma \ref{lem:ballinduct}.
For any $n > 0$, for any $x,x' \in V_n$, we have
$$
\sum_{i=0}^{n-1} \left| \log |Df(f^i(x))| - \log |Df(f^i(x'))|\right| \leq (1/2) \log 2 + \sum_{j =n(y)}^\infty e^{-n\delta} < \log 2.$$
Define $V_y$ as the set of $y' \in Y$ such that $\Pi F^{-n}(y') \in V_n$ for all $n \geq 0$. 

Then set $\alpha(y)$ equal to the radius of $V$. If $n(y) = 0$ set $\alpha(y) := 0$.  Then $y \mapsto \alpha(y)$ is a measurable function, positive almost everywhere. 

If $n(y)>0$, set
$$\rho(y) := \inf\lbrace t>2 : 
t^{-1}e^{n(\chi - \delta)} \leq |Df^n(\Pi F^{-n}y')| \leq te^{n(\chi + \delta)} \mbox{ for all } n\geq 0, y' \in V_y\rbrace.
 $$
   If $n(y)=0$, set $\rho(y) := 0$. Then $\rho(y)$ is defined everywhere   
 (for example by Lemma \ref{lem:sea} and the distortion bound) and is measurable. 
\eprf

\section{Good inverse branches}
In this section we define and show the existence of regularly returning cylinders in the natural extension. 
We start off with a fact which
 we borrow from \cite{Ledrappier:AbsCnsInterval} and a couple of lemmas.
\begin{Fact}{\cite{Ledrappier:AbsCnsInterval}}\label{Fact:messums}
Let $m$ be a finite measure on an interval $I$. For Lebesgue almost every $w$ in $I$, for all $\alpha < 1$, 
$$
\limsup_{n \to \infty} \frac{1}{n} \log m([w-\alpha^n, w + \alpha^n]) < 0.
$$
\end{Fact}
\begin{lem} \label{lem:gammasubexp}
Let $f : \cbar \to \cbar$ be a rational map with an invariant probability measure $\mu$ and natural extension $(Y, F, \mu)$. For each  $x \in \cbar$, for Lebesgue almost every $\gamma \in [0,|\cbar|]$, for almost all $y \in Y$
$$
\liminf_{n \to \infty} \frac{1}{n} \log \dist(\Pi F^{-n}y, \partial B(x,\gamma)) = 0;
$$
in words, the distance of $\Pi F^{-n}y$ to $\partial B(x,\gamma)$ decreases at fastest subexponentially.
\end{lem}
\beginpf
Recall that we denote by $\sigma(\cdot, \cdot)$ the spherical metric.
Define the measure $m$ on $[0, |\cbar|]$ by $m(C) := \omu( \{y : \sigma(\Pi y, x) \in C\})$. By invariance, for all integers $n$, 
$$
\omu( \{y : \sigma(\Pi y, x) \in C\})
 = \omu( \{y : \sigma(\Pi F^{-n} y, x) \in C\}).$$
 By Fact \ref{Fact:messums}, for all $0< \alpha < 1$, for almost all $\gamma$, 
$$
\limsup_{n \to \infty} \frac{1}{n} \log \omu( \{y : \sigma(\Pi F^{-n}y,x) \in [\gamma-\alpha^n, \gamma + \alpha^n]\}) < 0,
$$
so $\sum_n \omu(\{y: \dist(\Pi F^{-n}y, \partial B(x,\gamma)) < \alpha^n\})$ is finite. Apply the Borel-Cantelli Lemma to conclude.
\eprf 
\begin{lem} \label{lem:fixM}
Let $M > 0$ be given and denote by $\Fix(M!)$ the fixed points of $f^{M!}$. There exists a continuous function $\theta$ such that $\theta(0) = 0$ and such that, for any set $W$,  
  $W \cap f^{k}(W) \ne \emptyset$ and $0 < k \leq M$ imply that 
\begin{equation}\label{eqn:xiupe}
\dist(W, \Fix(M!)) \leq \theta(|W|).
\end{equation}
\end{lem}
\beginpf
Let $\varepsilon >0$ and let $C$ denote the complement of the set $B(\Fix(M!), \varepsilon)$. Each point $x$ in $C$ satisfies $\sigma(x,f^k(x)) > 0$ for each $1 \leq k \leq M$.
 Since $C$ is compact, there exists a $\delta < \varepsilon$ such that for these $k$ and for all $x \in C$, $\sigma(x,f^k(x)) > \delta$. Similarly, there exists a $\nu >0$ such that, if $W$ satisfies $|W| < \nu$, then  $f^k(|W|) < \delta/2$ for each $k, 0 \leq k \leq M$. Suppose $|W| < \nu$ and $f^{k}(W) \cap W \ne \emptyset$ for some $k$, $1 \leq k \leq M$. Then there is a point $x \in W$ such that $\sigma(x,f^k(x)) < \delta$. In particular, $x \in B(\Fix(M!), \varepsilon)$ and $\dist(W, \Fix(M!)) <  \varepsilon$.  
\eprf

\begin{dfn} \label{dfn:regretcyl} 
A \emph{regularly returning cylinder for $(Y, F, \omu)$ (and $\eta > 0$)} is a subset $A$ of $Y$ with the following properties:
\begin{itemize}
\item
$\omu(A) > 0$; 
\item
the base of $A$, $\Pi A$, is disjoint from the critical set and is a simply connected open set $U$ with smooth boundary;
\item there is a simply connected set $W \supset U$ such that $\dist(U, \partial W) > 0$ and 
for each $y$ in $A$, there are sets $U_y, W_y$ with  $ y \in U_y \subset W_y \subset Y$ such that  for all $n \geq 0$, $f^n : \Pi F^{-n} W_y \to W$ is a diffeomorphism, and $\Pi U_y = U$;
\item for all $y$ in $A$, for all $y' \in W_y$,  one has
$$
 \sum_{i=1}^\infty \left|\log |Df(\Pi F^{-i}y'))| -\log |Df(\Pi F^{-i}y)|\right|    < \log 2;
$$ 
\item there exists a constant $\rho_0 > 1$ such that, for all $y \in A$,
$$
\rho_0^{-1} e^{n(\chi_\mu - \eta)} \leq |Df^n(\Pi F^{-n}y)| \leq \rho_0 e^{n(\chi_\mu + \eta)}.$$
\item  for any $n, n' \geq 0$ and $y,y' \in A$, $\Pi F^{-n} U_y$ and $\Pi F^{-n'}U_{y'}$ are either nested or disjoint and 
  $F^{-n} U_y$ and $F^{-n'}U_{y'}$ are either nested or disjoint;
\item  for all $n \geq n' \geq 0$ and $y, y' \in A$, $\Pi F^{-n} W_y  \cap \partial \Pi F^{-n'}U_{y'} = \emptyset$.
\end{itemize}
\end{dfn}
\begin{thm} \label{thm:regretcyl}
Let $f : \cbar \to \cbar$ be a rational map, $\mu \in \scrM^+(f)$  and  $(Y,F,\omu)$  the natural extension. Then there exists a regularly returning cylinder for $(Y,F,\omu)$.
\end{thm}
\beginpf 
%
Let $\alpha(y)$ be given by 
 Proposition \ref{prop:complexNeilUnstableManif}. 
There exists  $\alpha_0 > 0$ such that, setting  
$$
A_0:= \{y \in Y : \alpha(y) > \alpha_0\},
$$
$A_0$
has positive measure. Hence there is an $x \in \cbar$ such that $A_1 := A_0\cap  \Pi^{-1}B(x, \frac{\alpha_0}{9})$ has positive measure. 

By  Lemma \ref{lem:gammasubexp}, there exists a $\gamma$ satisfying $\frac{\alpha_0}{9} < \gamma < \frac{\alpha_0}{8}$  such that, for almost every $y \in Y$, $\dist(\Pi F^{-n}y, \partial B(x, \gamma))$ decreases sub-exponentially. Note that this implies $\mu(\partial B(x,\gamma)) = 0$.

Set $U' := B(x,\gamma)$, then $\overline{U'} \subset B(x, \frac{\alpha_0}{8}) \subset B(x, \frac{3\alpha_0}{4}) \subset B(\Pi y, \alpha(y))$ for all $y$ in $A_1$. 
Let $\eta$ satisfy $0 < \eta < \chi_\mu$. Let $V_y$ and $\rho$ be given by  Proposition \ref{prop:complexNeilUnstableManif}, so for all $y' \in V_y$,
$$
|Df^n(\Pi F^{-n}y')| \geq \rho(y) e^{n(\chi_\mu - \eta)}
$$
for all $n \geq 1$. For $t > 0$, set 
$$
A_t := 
  \lbrace y \in A_1 : \mbox{ for all } n > 0, \dist(\Pi F^{-n}y,\partial U') \geq te^{-n\frac{\chi_\mu - \eta}{2}} \hfill 
$$
$$
\hfill 
\mbox{ and } |Df^n_{|\Pi F^{-n}V_y}| \geq t e^{n(\chi_\mu -\eta)} \rbrace.
$$ 
Since $A_1 = \bigcup_{t>0} A_t$ up to a set of zero measure, 
 there is a $t_0 > 0$ such that $A_2 := A_{t_0}$ has positive measure. For some $M >0$ sufficiently large, for all $y \in A_2$ and $n \geq M$, 
$$
|\Pi F^{-n}V_y| < |\cbar| \frac{1}{t_0} e^{-n(\chi_\mu - \eta)} < t_0 e^{-n\frac{\chi_\mu - \eta}{2}}.
$$
 Thus for all $y \in A_2$, for all $n \geq M$, 
\begin{equation}\label{eqn:xiuprime}
\Pi F^{-n}V_y \cap \partial U' = \emptyset.
\end{equation}
We are halfway there. We have shown that if you go \emph{far enough back} along an inverse branch you will no longer intersect $\partial U'$. 

%

 Pulling back (\ref{eqn:xiuprime}), if $n \geq 0$ and there exists a $k > 0$ such that 
\begin{equation}\label{eqn:xiupe2}
 \Pi F^{-n-k}V_y \cap \partial \Pi F^{-n}U_y'  \ne \emptyset,
\end{equation}
then $k \leq M$.

For each $y$ in $A_2$, let $\{n_j(y)\}$ be the strictly increasing collection of $n \geq 0$ for which there exists a $k > 0$ satisfying
$$
 \Pi F^{-n-k}V_y \cap \partial \Pi F^{-n}U'_y  \ne \emptyset.
$$
 Note that $\{n_j(y)\}$ may be empty or finite. It is immediate that if $n_i(y), n_{i+1}(y) \in \{n_j(y)\}$, then $n_{i+1}(y) - n_i(y) = 1$, and we showed that for each such $n$ the minimal $k$ is $\leq M$.

Suppose  first there is a positive measure subset $C$ of $A_2$ such that for each $y$ in $C$ the set $\{n_j(y)\}$ is infinite (and thus equal to $\N$). Then (\ref{eqn:xiupe}) gives that
$$
\dist(\Pi F^{-n} V_y, \Fix(M!)) \leq \theta(|\Pi F^{-n} V_y|).
$$
But $|\Pi F^{-n} V_y|$ is exponentially decreasing, so $\dist\left(\Pi F^{-n} V_y, \Fix(M!)\right)$ tends to zero (and this uniformly in $y\in A_2$), by Lemma \ref{lem:fixM}. Thus $\mu\left(B(\Fix(M!), \varepsilon)\right) \geq \omu(F^{-n}C) = \omu(C)$ for all $\varepsilon > 0$. This implies that $\mu$ is atomic and supported on a repelling periodic orbit. In this case the theorem is easy. 

Otherwise,
 for almost every $y$ in $A_2$, the set $\{n_j(y)\}$ is finite; let $n(y)$ denote its maximum (defined as zero if the set is empty). 
Then there exists an $N > 0$ such that the set 
$$
A_3 := \{y \in A_3 : n(y) < N\}
$$
is of positive measure. Summarising, we have shown that, for all $y$ in $A_3$,
\begin{equation}\label{eqn:xeijeiupe}
 \Pi F^{-n-k}V_y \cap \partial \Pi F^{-n}U_y'  = \emptyset,
\end{equation}
for all $n \geq N$ and all $k > 0$.


Let $U$ be a connected component of $f^{-N} (U')$ such that 
$$
A_4 := F^{-N}A_3 \cap \Pi^{-1} U
$$
has positive measure (thus $U = \Pi F^{-N} U'_{F^Ny}$ for all $y \in A_4$).
Similarly let $W$ be the connected component of $f^{-N}(B(x,\frac{3\alpha_0}{4}))$ containing $U$. For each $y \in A_4$, set $W_y := (F^{-N} V_{F^Ny}) \cap \Pi^{-1}W$ and $U_y := W_y \cap \Pi^{-1} U = F^{-N}U'_{F^Ny}$.

Recall that 
$$
f^N(U) = U' \subset B(x,\alpha_0/8) \subset B(x,3\alpha_0/4) =f^N(W).$$
Since distortion of $f^N$ restricted to $W$ is bounded by $\log 2$, it follows that $\dist(U, \partial W) > |U| > 0$. 

By (\ref{eqn:xeijeiupe}), for all $n \geq 0$, and all $y \in A_4$,
$$
\partial U \cap \Pi F^{-n} W_y = \partial U \cap \Pi F^{-n} U_y = \emptyset.
$$
This is nearly a ``regularly returning'' property: It implies that for any $n\geq n' \geq 0$ and $y,y' \in A_4$, $\Pi F^{-n} W_y$ and $\Pi F^{-n'}U_{y'}$ are either nested or disjoint. Set 
$$
A := \bigcup_{y \in A_4} U_y
$$
and the theorem follows.
%
%
\eprf

%

\section{Partitions}\label{sec:partitions}
First we show that there is a nice partition of $Y$. Let $A$ be a regularly returning cylinder for $(Y,F,\omu)$. For $y \in Y$ let $e(y)$ be the first entry time of $y$ into $A$, that is, $e(y) := \inf\{n \geq 0: F^ny \in A\}$. This is defined (finite) almost everywhere by ergodicity, and one can set 
$$
\xi(y) := F^{-e(y)}U_{F^{e(y)}y}.$$
If $y$ does not enter $A$, set $\xi(y):= \emptyset$. By the regularly returning property, $e(y) = e(y')$ and $\xi(y) = \xi(y')$ for all $y' \in \xi(y)$. Thus one can define a partition $\xi$ of $Y$ where the element of $\xi$ containing $y$ is  $\xi(y)$ (up to a set of measure zero --- alternatively one can add back the complement, a set of measure zero).

\begin{prop} \label{prop:nicenice}
Given any  regularly returning cylinder $A$ for $(Y,F,\omu)$, there exists a measurable partition $\xi$ of $Y$ such that, if $\xi(y)$ denotes the element of $\xi$ containing $y$ and $e(y)$ the first entry time of $y$ to $A$, then 
$$
F^{e(y)}\xi(y) = U_{F^{e(y)}y}.$$
For any $y, y'$ in $Y$, $\Pi \xi(y)$ and $\Pi \xi(y')$ are either nested or disjoint, and $\xi(y)$ and $\xi(y')$ are either nested or disjoint.
\end{prop}
\beginpf Immediate.
\eprf

Now we can define a first partition $\zeta$ of $\cbar$ (modulo a set  of measure zero). Let 
$$
\zeta(z) := \bigcup_{\Pi y = z} \Pi \xi(y).$$
This is a countable partition of $\cbar$ with the property  that each element gets mapped by an iterate of $f$ diffeomorphically (and with bounded distortion) to $U = \Pi A \in \zeta$. Also, any element of $Q \in \zeta\setminus \{U\}$ has the property that $f(Q)$ is contained in an element of $\zeta$. 
We would now like to construct a finite generating partition with good topological properties. 

\begin{prop}\label{prop:nicebelow}
Let $f$ be a rational map of degree $d$ and let $\mu \in \scrM^+(f)$.
Let $A$ be a regularly returning cylinder for $(Y,F,\omu)$ with base $U$. Then there exists a measurable partition $\scrP$ with the following properties: 
\begin{itemize}
\item $\scrP$ has at most $d+1$ elements $P_0, \ldots, P_d$;
\item $f$ is injective on each partition element; 
\item one partition element $P_0 \in \scrP$ is the base $U$; 
\item $f^{-1}(f(U)) \cap P_d = \emptyset$;
\item each connected component of a partition element is mapped by an iterate of $f$ diffeomorphically to $U$, and moreover is the projection of an element of the partition $\xi$ of Proposition \ref{prop:nicenice}; 
\item the projection of each element of $\xi$ is contained in exactly one element of $\scrP$;
\item the partition is generating.
\end{itemize}
\end{prop}



\beginpf We shall define the partition elements $P_0, \ldots, P_d$ of the partition $\scrP$ by distributing the elements of $\zeta$ to them. Firstly set $P_0 := U$. 
To each $Q \in \zeta$ there is a unique $n_Q \geq0$ such that $f^{n_Q}Q = U$. Let $\zeta_n := \{Q \in \zeta : n_Q = n\}$. 
For all $n \geq 1$, for all $Q \in \zeta_n$, one has $f(Q) \subset Q' \in \zeta_m$ for some $Q',m$ satisfying $0\leq m < n$. 

Suppose $P_0, P_1^n, P_2^n, \ldots, P_d^n$ is a partition of $\bigcup_{m< n}\bigcup_{Q\in \zeta_m} Q$ such that 
$f$ is injective on each partition element,
 initialising $P_i^0 := \emptyset$ for $1 \leq i \leq d$.

\begin{itemize}
\item
Consider an element $Q$ of $\zeta_n$ and let $\zeta_Q$ be collection of elements $Q' \ni \bigcup_{j\leq n} \zeta_j$ such that $f(Q') \cap f(Q) \ne \emptyset$. It follows that for $Q' \in \zeta_Q$, $f(Q') \supset f(Q)$.
Since $f$ is  $d$-to-one, there are at most $d$ elements in $\zeta_Q$. The elements of $\zeta_Q \setminus \zeta_n$ are already subsets of  distinct partition elements $P_i$. We distribute the remaining ones (contained in $\zeta_Q \cap \zeta_n$), so that each element of $\zeta_Q$ is a subset of a different $P_i^{n+1}$: 

\item
For $1 \leq i \leq d$ let $P'_i := P_i^n$.
\item
Order $\zeta_Q \cap \zeta_n =: Q_1, \ldots Q_k$. Let $i_1 := \inf\{i \geq 1: f(Q_1) \cap f(P_i) = \emptyset\}$.  Let $i_2 := \inf\{i > i_1 : f(Q_2) \cap f(P_i) = \emptyset\}$. Continue inductively and replace $P'_{i_j}$ by $P'_{i_j} \cup Q_j$ for $j = 1, \ldots, k$. It follows from the construction that $f$ remains injective on each $P'_i$ and $U \cap P'_d = \emptyset$. 
\item
Repeat until there are no further elements in $\zeta_n$, that is, until for all $Q \in \zeta_n$ there exists an $i$ such that $Q \subset P'_i$. 
\item
Set $P_i^{n+1} := P'_i$ for $i = 1, \ldots, d$. 
\end{itemize}

We proceed inductively, exhausting all of $\zeta$, setting
$P_i:= \bigcup_{n\geq0} P_i^n$ for $i = 1, \ldots d$, and $\scrP := \{P_0, P_1, \ldots, P_d\}$. 

It remains to show that the partition is generating, that is that on a set of full measure the partitions $\scrP_n := \bigvee_{i=0}^n f^{-i} \scrP$ tends to the partition into points. Let $X$ be a set of full measure such that for all $x \in X$ there is a $y \in \Pi^{-1}x$ such that $y$ enters $A$ infinitely often (this set $X$ exists by ergodicity and the positivity of the measure of $A$).

Suppose $x,x' \in X$ are in the same element of $\scrP_n$ for all $n \geq 0$. There are a $y \in \Pi^{-1}x$ and an infinite sequence $n_j \nearrow \infty$ such that $F^{n_j}y \in A$. Thus there is an infinite sequence $U_{j}$ such that $f^{n_j}: U_{j} \to U \in \scrP$ is a diffeomorphism. For each $0 \leq i \leq n_j$, $f^i(U_{j})$ is contained in an element of $\zeta$ and thus in an element of $\scrP$.
Now $f$ is injective on each partition element, so 
$$
U_j \in \scrP_{n_j}.$$ 
We also have by Theorem \ref{thm:regretcyl} that $|U_j| < \rho_0 e^{-n_j(\chi -\eta)} |U|$, so $|U_j| \searrow 0$ and $x = x'$.
\eprf

\begin{lem} \label{lem:belowproperties}
Let  $A$ be a regularly returning cylinder with base $U$.
Let $\scrP$ be the partition of $\cbar$ given by Proposition \ref{prop:nicebelow}.  Then for each $x$ in a  set $X$ of full measure, there exists a strictly monotone sequence $n_j \nearrow \infty$ and corresponding sets $U_j \ni x$ such that 
\begin{itemize}
\item $\lim_{j \to \infty} \frac{n_j}{j} > 0$, so
$$ 
\lim_{n \to \infty} \frac{n_j}{n_{j+1}} = 1;$$
\item  $f^{n_j}: U_j \to U$ is a diffeomorphism with uniformly bounded distortion;
\item $\lim \frac{1}{n} \log|Df^{n_j}(x)| = \chi$;
\item $\dist(f^{n_j}(x), \partial U) \geq \varepsilon $ for some $\varepsilon > 0$;
\item $U_j$ is an element of the partition $\scrP_{n_j} = \bigvee_{i=0}^{n_j} f^{-i}\scrP$.
\end{itemize}
\end{lem}
\beginpf
We have essentially proven this already. Let $B$ be a subset of $A$ such that $\omu(B) > 0$ and $\dist(\Pi B, \partial U) \geq \varepsilon$ for some $\varepsilon > 0$. Almost every point $x$ is the projection of a point $y \in Y$ which enters $B$, by ergodicity, at successive times $n_j$ with $\lim_{j \to \infty} j/n_j = \omu(B)$. Then the element of the partition $\scrP_{n_j}$ containing $x$ is $\Pi F^{-n_j} U_{F^{n_j}y}$. The uniform bound on the distortion comes from  the definition of regularly returning cylinders. Since $F^{n_j}y \in B$, one has $\dist(f^{n_j}(x), \partial U) \geq \varepsilon$. 
\eprf

\section{Dimensions of measures}\label{sec:dvl}
In this section we give short proofs of two results concerning the dimension of a measure, making use of the partitions of Section \ref{sec:partitions}.
The Hausdorff dimension of a measure, denoted $\HD(\mu)$, is defined as the infimum of Hausdorff dimensions of sets of full measure.  
In dimension one, a theorem proving that 
$$
\HD(\mu) = \frac{h_\mu}{\chi_\mu}
$$
for an ergodic invariant probability measure with positive Lyapuonv exponent for some class of maps is called a \emph{Dynamical Volume Lemma}. Our proof of the following well-known result is shorter than usual as we have at our disposal a finite generating partition with good properties. 
\begin{prop}[\cite{Ledrappier:AbsCnsComplex}]  \label{prop:ratdvl}
Let $f:\cbar \to \cbar$ be a rational map and $\mu \in \scrM^+(f)$.  Then for almost all $x$ with respect to $\mu$,
$$
\lim_{r \to 0} \frac{\log \mu(B(x,r))}{\log r} = \frac{h_\mu}{\chi_\mu}.
$$
In particular, 
$$
\mbox{HD}(\mu) = \frac{h_\mu}{\chi_\mu}.
$$
\end{prop}
\beginpf
We write $\chi$ for $\chi_\mu$. 
Let $(Y, F, \overline{\mu})$ be the natural extension. Let $\eta$ be a small positive constant. Let $A$ be a regularly returning cylinder given by Theorem \ref{thm:regretcyl}, $\xi$ a corresponding partition of $Y$ given by Proposition \ref{prop:nicenice} and $\scrP = \{P_0 = U = \Pi A, P_1,\ldots, P_d\}$ a corresponding finite generating partition of $\cbar$ given by Proposition \ref{prop:nicebelow}.
The result essentially follows from Lemma \ref{lem:belowproperties}: let $X$ be a set of full measure given by this lemma.

Let $\scrP_n := \bigvee_{i=0}^n f^{-i} \scrP$ and $\scrP_n(x)$ the element of $\scrP_n$ containing $x$.
By the Shannon-McMillan-Breiman Theorem (\cite{Parry:Book} p.39), there is  a set $X' \subset X$ of full measure such that 
$$
h_\mu = \lim_{n \to \infty} \frac{-1}{n}\log \mu(\scrP_n(x))$$
 for all $x \in X'$. Now fix $x \in X'$ and let $n_j$ and $U_j$ be given by Lemma \ref{lem:belowproperties}. Then $U_j = \scrP_{n_j}(x)$.
Thus
$$
h_\mu = \lim_{j \to \infty} \frac{-1}{n_j}\log \mu(U_j).$$

Now the uniform bound on the distortion together with $\dist(f^{n_j}(x), \partial U) \geq \varepsilon >0$ implies that there are annuli $A_j \supset \partial U_j$ centred on $x$ such that the moduli $A_j$ are uniformly bounded. Let $r_j$ be the inner radius and $R_j$ the outer, and $K$ such that $R_j \leq K r_j$. 
Since 
 $\lim_{j \to \infty} \frac{n_j}{n_{j+1}} = 1$, we have
$$
\lim_{j \to \infty} \frac{1}{n_j}\log r_{j+1} = 
\lim_{j \to \infty} \frac{1}{n_{j+1}}\log R_j = 
-\chi.$$
Thus, if $r_j \geq r \geq r_{j+1}$,
$$
\frac{\log\mu(B(x,r))}{\log r} \geq \frac{\log\mu(U_j)}{\log r_{j+1}} 
=  \frac{-1}{n_j}\log{\mu(U_j)} \frac{-n_j }{\log r_{j+1}} 
$$
and the right-hand side tends to $\frac{h_\mu}{\chi}$ as $j\to \infty$.
If $R_j \geq r \geq R_{j+1}$, 
$$
\frac{\log\mu(B(x,r))}{\log r} \leq \frac{\log\mu(U_{j+1})}{\log R_{j}} 
= \frac{-1}{n_{j+1}}\log{\mu(U_{j+1})} \frac{-n_{j+1}}{\log R_j} 
$$
and the right-hand side tends to $\frac{h_\mu}{\chi}$ as $j\to \infty$.

As $r \to 0$ one has $j \to \infty$ so
we conclude that 
$$
\lim_{r \to 0} 
\frac{\log\mu(B(x,r))}{\log r} 
= \frac{h_\mu}{\chi}$$
as required.
\eprf

The following proposition also has a surprisingly short proof using the same techniques.
\begin{prop}\label{prop:complexHDmu}
Let $f$ be a rational map of the Riemann sphere and $m$ a non-exceptional $(\phi,t)$-conformal measure for $f$. Suppose  $\mu \in \scrM^+(f)$  is absolutely continuous with respect to $m$. 

Then $\mu$ and $m$ are equivalent and HD$(\mu) = t + \frac{\int \phi d\mu}{\chi_\mu}$.
\end{prop}
\beginpf
First let us show equivalence of the measures. Suppose there exists a set $B$ such that $m(B) > 0$ but $\mu(B) = 0$. Then 
$$X := \cbar \setminus \bigcup_{i\geq 0} f^{-i}(B)
$$
is a forward-invariant set of full $\mu$-measure disjoint from $B$. Combining the density point theorem with Lemma \ref{lem:belowproperties}, for $\mu$ almost every point $x \in X$, there exist $U_j$ such that
$$
\lim_{j \to \infty} \frac{m(X \cap U_j)}{m(U_j)} = 1
$$
and $f^{n_j} : U_j \ni x : \to U$ has uniformly bounded distortion. Thus $m(X \cap U) = m(U)$. Then the eventually onto property implies $m(\cbar \cap X) = m(\cbar)$ (we use here that $m$ is non-exceptional). But then $m(B) = 0$, contradiction. 

Now it suffices to show the dimension estimate for $m$ instead of $\mu$. Let $x$ be a typical point (for $m$ or for $\mu$, it is the same thing) and let $U_j$, $r_j$ and $R_j$ be as per the proof of the preceding proposition. 
As before, we get
$$
\lim_{j\to \infty} \frac{\log m(U_j)}{\log r_{j+1}} \leq \lim_{r\to 0} \frac{\log m(B(x,r))}{\log r} \leq \lim_{j\to \infty} \frac{\log m(U_{j+1})}{\log R_j},
$$
provided the limits exist. 

It now suffices to show that $\lim_{j \to \infty} (-1/n_j)\log m(U_j) = \chi_\mu t + \int \phi d\mu$ for almost every $x$. But conformality implies that 
$$
m(U) = \int_{U_j} e^{S_{n_j}\psi(x')} dm(x'),$$ 
where $S_n\psi := \sum_{i=0}^{n-1} (t\log|Df| +\phi) \circ f^i$. By Lemma \ref{lem:cvgphi} the integrand  has uniformly bounded variation on the sets $U_j$.  Thus 
$$
m(U) \approx m(U_j) e^{S_{n_j} \psi(x)}
 $$
and so
$$
\lim_{j\to\infty} (-1/n_j) \log m(U_j) = \lim_{j\to\infty}(1/n_j) S_{n_j}\psi(x).$$
By ergodicity, the right-hand term tends to $t\chi_\mu + \int \phi d\mu$. 
\eprf

\section{Entropy}
The goal of this section is to show that $h_\mu = H(F^{-1}\xi|\xi)$ for the partition $\xi$ of the natural extension given by Proposition \ref{prop:nicenice}.  
This partition is not (or does not seem to be) the future of some finite generating partition (which would imply the result), so there is something to be proved here. 
%
\begin{prop}\label{prop:niceentropy}
Given a rational map $f$, $\mu \in \scrM^+(f)$ and natural extension $(Y,F,\omu)$, let $A$ be a regularly returning cylinder 
and $\xi$ the corresponding partition of $Y$ given by Proposition \ref{prop:nicenice}. Then 
$$
h_\mu = h_\omu = H(F^{-1}\xi | \xi).$$
\end{prop}
\beginpf
We refer to \cite{Rohlin:Exact} and \cite{Walters:ErgodicBook} for the notions from ergodic theory used in this section, and to \cite{Walters:ErgodicBook} for notation. By section 3.3 of \cite{Rohlin:Exact}, $h_\mu = h_\omu$.
Let $\scrP$ be the finite generating partition given by Proposition \ref{prop:nicebelow} and $\Pi^{-1}\scrP$ the induced finite generating partition of $Y$, so $H(F,\Pi^{-1}\scrP) = h_\mu$. Set
$$
\zeta := \bigvee_{i=0}^{\infty} F^i \Pi^{-1}\scrP.$$
Then $H(F^{-1}\zeta | \zeta) = h_\mu$. We want to compare $H(f^{-1}\xi|\xi)$ with $H(f^{-1}\zeta|\zeta)$.
Associated to each partition is a canonical system of conditional probability measures, also called the Rohlin decomposition of $\omu$ with respect to that partition. 
Let $p_\xi(y,\cdot)$ be the decomposition of $\omu$ with respect to $\xi$, so for (almost) all $y$, $p_\xi(y, \cdot)$ is a probability measure on  $\xi(y)$ and for all $X \subset Y$,
$$
\omu(X) = \int_Y p_\xi(y,X\cap \xi(y)) d\omu.$$
Let $p_\zeta(y,\cdot)$ be the Rohlin decomposition of $\omu$ with respect to $\zeta$. 
For  $y \notin A$, $\xi(y) = F^{-1}(\xi(Fy))$, so $p_\xi(y,F^{-1}(\xi(Fy))) = 1$. Thus
$$
H(F^{-1}\xi|\xi) = \int_Y \log p_\xi(y, F^{-1}(\xi(Fy))) = \int_A \log p_\xi(y, F^{-1}(\xi(Fy))).$$
On $A$, $\xi(y) = \zeta(y) = U_y$, the leaf of $A$ containing $y$. Off $A$, $\zeta(y)$ is a union of elements of $\xi$. 

\begin{eqnarray*}
H(F^{-1} \zeta | \zeta) 
&=& \int_Y \log p_\zeta(y, F^{-1}(\zeta(Fy)) d\omu \\
&=&
\int_A \sum_{i=0}^{r(y) -1} \log p_\zeta(F^iy, F^{-1}(\zeta(F^{i+1}y))) d\omu\\
&=& \int_A \log \prod_{i=0}^{r(y) -1} p_\zeta(F^iy, F^{-1}(\zeta(F^{i+1}y))) d\omu\\
&=& \int_A \log p_\zeta(y, F^{-r(y)}(\zeta(F^{r(y)}y))) d\omu\\
&=& H(F^{-1}\xi|\xi).
\end{eqnarray*}
\eprf

\section{Absolute continuity}
 Let $\phi : \cbar \to  \arr$ be a H\"older continuous function and $C,\epsilon >0$ such that $|\phi(x) - \phi(x')| \leq C\sigma(x,x')^\epsilon$ for all $x,x' \in \cbar$. Since $\cbar$ is compact, $\phi$ is necessarily bounded. 
If $\psi : \cbar \to \arr$ is any function, we write 
$$
S_n\psi(x) := \sum_{i=0}^{n-1}\psi \circ f^i(x)$$
 and, for $y \in Y$, 
$$
(S_{-n}\psi)(y) := \sum_{i=1}^{n} \psi \circ \Pi \circ F^{-i} (y).$$
\begin{lem}\label{lem:cvgphi}
Let $A$ be a regularly returning cylinder for $(Y,F,\omu)$ given by Theorem \ref{thm:regretcyl} and $\xi$ the corresponding partition given by Proposition \ref{prop:nicenice}. Let $t \in \arr$ and set $\psi := \phi +  t \log|Df|$.
There exists a constant $C >0$ such that for almost every $y \in Y$, for all $y' \in  \xi(y)$ and for all $n \geq 0$,
$$
(S_{-n} \psi)(y) - (S_{-n} \psi)(y')
\leq \sum_{i=1}^\infty \left|\psi \Pi F^{-i}y - \psi \Pi  F^{-i}y'\right| \leq C, $$ 
and so, taking exponentials and defining $\Phi(y,\cdot) : \xi(y) \to \arr$ by
$$
\Phi(y,y') := \lim_{n \to \infty} e^{(S_{-n}\psi)(y') - (S_{-n}\psi)(y)},
$$
$\Phi(y,\cdot)$ is well-defined and uniformly bounded away from zero and infinity. 
\end{lem}
\beginpf
Note that the first inequality is immediate, by definition, and that the second inequality implies that the series 
$\sum_{i=1}^\infty (\psi \Pi F^{-i}y - \psi \Pi F^{-i}y')$ 
is absolutely convergent,  so $\lim_{n\to \infty} 
(S_{-n} \psi)(y) - (S_{-n} \psi)(y')$ exists and is finite. In particular, we only need to show the second inequality.

By  the definition of regularly returning cylinders, 
$$
\sum_{i=1}^\infty \left|t\log|Df(\Pi F^{-i}y)| - t\log|Df(\Pi F^{-i}y')|\right|  \leq |t|\log 2.$$
Since $|\psi(x) - \psi(x')|\leq |\phi(x)-\phi(x')| + t|\log |Df(x)| - \log|Df(x')||$, 
it only remains to treat 
$$
\sum_{i=1}^\infty |\phi (\Pi F^{-i}y) - \phi (\Pi F^{-i}y')| 
$$
But, using H\"older continuity and the definition of regularly returning cylinders, 
$$
|\phi (\Pi F^{-i}y) - \phi (\Pi F^{-i}y')| \leq C \sigma(\Pi F^{-i}y,  \Pi F^{-i}y')^\epsilon$$
$$
\leq C  \rho_0^\epsilon  e^{-i(\chi-\eta)\epsilon} |U|^\epsilon$$
so
$$
 \sum_{i= 1}^\infty |\phi (\Pi F^{-i}y) - \phi (\Pi F^{-i}y')| 
\leq
C \rho_0^\epsilon \sum_{i\geq1} e^{-i(\chi-\eta)\epsilon} |U|^\epsilon 
 \leq M  < \infty
$$
for some constant $M$ independent of $y$, as required.
\eprf

The Rohlin decomposition $p(y,\cdot)$ for the measure $\overline{\mu}$ with respect to the partition $\xi$ is a conditional probability measure (the conditional expectation) on each partition element of $Y$ such that, for any measurable set $B \subset Y$ one has
$$
\overline{\mu}(B) = \int_Y p(y,B) d\overline{\mu} = \int_Y p(y, B\cap\xi(y)) d\overline{\mu}.
$$
For all $n > 0$, by Proposition \ref{prop:niceentropy},
\begin{equation} \label{eqn:ratnentropy}
nh_\mu = H(F^{-n}\xi | \xi) = -\int \log p(y, [F^{-n}\xi](y)) d\overline{\mu},
\end{equation}
where $[F^{-n}\xi](y)$ denotes the element of the partition $F^{-n}\xi$ containing $y$. 

We shall now start to study $(\phi,t)$-conformal measures.
In the following proposition, $q(y,\cdot)$ is an educated guess as to the Rohlin decomposition: it could be derived directly if we knew in advance that $\mu(B) = \lim_{n \to \infty} m(f^{-n}B)$ for all measurable sets $B$. 
\begin{prop}\label{prop:jimothy}
    Let $f:\cbar \to \cbar$ be a rational map and $\mu \in \scrM^+(f)$. Let $(Y,F,\omu)$ be the natural extension and $A$ be a regularly returning cylinder and $\xi$ the corresponding partition given by Proposition \ref{prop:nicenice}. 

Suppose that there is a non-exceptional $(\phi,t)$-conformal measure $m$ such that $h_\mu = t\chi + \int \phi d\mu$. 
Let $\psi$ and $\Phi$ be defined as per Lemma \ref{lem:cvgphi}. 

The Rohlin decomposition of $\omu$ with respect to the partition $\xi$ is given by $q(y,\cdot)$, where 
\begin{equation}\label{eqn:ratrohlin}
q(y,dz) := \frac{ \Phi(y,z) dm_{\xi(y)}(z)}{\int_{\xi(y)} \Phi(y,y') dm_{\xi(y)}(y')},
\end{equation}
where $m_{\xi(y)}$ is the pullback by $\Pi_{|\xi(y)}$ of the conformal measure $m$ restricted to $\Pi\xi(y)$.
\end{prop}
\beginpf 
 The proof of this proposition is very similar to that of proposition 3.6 of \cite{Ledrappier:AbsCnsInterval}; we include it for the reader's convenience.  

Firstly, note that since $\chi_\mu > 0$, $\mu$ is supported on the Julia set, and $m$ is non-exceptional, so $m$ is positive on the $\Pi \xi(y)$ for $\omu$-almost every $y$.
Since $\Phi$ is uniformly bounded away from zero, $q(y,\cdot)$ is well-defined almost everywhere. 

It follows from (\ref{eqn:ratrohlin}) that $q(y,\cdot)$ is a probability measure on $\xi(y)$, and that if $\xi(y_1) = \xi(y_2)$ then $q(y_1,\cdot) = q(y_2, \cdot)$ by multiplying above and below by $\Phi(y_1,y_2)$.
Set $k(y) := \int_{\xi(y)} \Phi(y,y') dm(y')$. Then, by definition,
\begin{eqnarray*}
q(y, [F^{-n}\xi](y)) &=& \frac{\int_{[F^{-n}\xi](y)} \Phi(y,y') dm(y')}{k(y)} \\
&=& \frac{1}{k(y)}
\int_{F^{-n}(\xi(F^n(y)))} \Phi(F^ny,F^ny') 
e^{S_n\psi(\Pi y')} e^{ - S_n\psi(\Pi y)}\\\
&=& \frac{k(F^ny)e^{-S_n\psi(y)}} {k(y)}.
\end{eqnarray*}
We used that the Jacobian for $f^n$ of the conformal measure $m$ is $e^{S_n \psi}$. Taking logs,
$$
\log q(y, [F^{-n}\xi](y)) = 
\log k(F^ny) - t\log |Df^n(\Pi y)| -S_n\phi(y) - \log k(y).
$$
Integrating, applying the classical lemma 3.5 of \cite{Ledrappier:AbsCnsInterval}, 
$$
-\int \log q(y, [F^{-n}\xi](y)) d\overline{\mu} = nt \chi_\mu + n \int \phi d\mu.
$$
Let $p(y, \cdot)$ be the Rohlin decomposition for the measure $\overline{\mu}$ with respect to the partition $\xi$. 
Recalling (\ref{eqn:ratnentropy}), we deduce
$$
\int \log q(y, [F^{-n}\xi](y)) d\overline{\mu} = \int \log p(y, [F^{-n}\xi](y)) d\overline{\mu},
$$
so
\begin{equation} \label{eqn:ratlogrohlin}
\int \log \frac{q(y, [F^{-n}\xi](y))}{p(y, [F^{-n}\xi](y))} d\overline{\mu}.
\end{equation}
Consider the measure $\overline{\nu}_n$ defined on the $\sigma$-algebra of $[F^{-n}\xi]$ measurable sets defined by 
$$
\overline{\nu}_n(B) = \int q(y, B) d \overline{\mu}.
$$
By (\ref{eqn:ratlogrohlin}), 
$$
\int \left(\log \frac{d\overline{\nu}_n}{d\overline{\mu}}\right) d\overline{\mu} = 0.
$$
Both measures are probabilities, so by concavity of logarithm, $\overline{\nu}_n = \overline{\mu}$, so $q(y, B) = p(y, B)$ if $B$ is $[F^{-n}\xi]$ measurable. As $\xi$ generates and $\xi < F^{-1}\xi$, one concludes that $q$ and $p$ coincide almost everywhere.
\eprf
\begin{cor}\label{cor:pesinac}
Under the same hypotheses, $\mu$ is absolutely continuous with respect to the $(\phi,t)$-conformal measure $m$.
\end{cor}
\begin{prop}\label{prop:acmarkov}
Let $f:\cbar \to \cbar$ be a rational map, $\mu \in \scrM^+(f)$ and  $(Y,F,\omu)$ the natural extension. 

Suppose $m$ is a non-exceptional $(\phi,t)$-conformal measure  and that $h_\mu = t\chi_\mu + \int \phi d\mu$. 
Then  the density of $\mu$ is bounded from below by a constant $\varepsilon >0$ on an open set $U$ 
and there exists an expanding induced Markov map for $(f,m)$ with integrable return time which generates $\mu$.
If $t \geq 0$ then the density of $\mu$ is bounded away from zero  almost everywhere with respect to $m$. 
\end{prop}
\beginpf
Let $A$ be a regularly returning cylinder with base $U$ and let $\xi$ be the corresponding partition.  
Let $\omu_A := \omu_{|A}$. Since $\Phi$, as defined in Lemma \ref{lem:cvgphi},  is uniformly bounded away from zero and infinity on $A$, the density of $\Pi_* \omu_A$ is bounded away from zero and infinity on $U$. In particular it is bounded from below by some constant $\varepsilon >0$. Since $\Pi_* \omu_A(C) \leq \omu(\Pi^{-n}C) = \mu(C)$, the density $\rho$ of $\mu$ is bounded from below by $\varepsilon >0$ on $U$. 

By the eventually onto property, there exists a finite $n$ such that $f^n(U) \supset \J$. By invariance of $\mu$ and conformality of $m$, for each $x \in U$,  
$$
\rho(f^n(x)) \geq 
\varepsilon \left(|Df^n(x)|^t e^{S_n \phi(x)} \right)^{-1}.$$
In particular, 
if $t\geq 0$ then the density is bounded away from zero everywhere on $\J$ and thus on $\cbar$.

Let $N>0$ be large enough that for all $n \geq N$, for all $y \in A$, $|Df^n(\Pi F^{-n}y)| > 2$ and for each $y \in A$ let $m(y)$ denote the $N^{\mathrm{th}}$ return time of $y$ to $A$. In particular $m(y) > N$ almost everywhere on $A$, and  by a simple extension of Kac's Theorem (\cite{Petersen:ErgodicBook}, theorem 2.4.6),   
$$
\int_A m(y) d\omu = N.
$$
Then $\{\Pi F^{-m(y)} \xi(F^{m(y)}) : y \in A\}$ is a cover of a subset $X \subset U$ of full measure by nested or disjoint topological balls which get mapped by an iterate of $f$ diffeomorphically onto $U$ with extensions that get mapped onto $V$. By the nested or disjoint property, there is a pairwise disjoint countable cover $\{U_i\}_{i\in \N}$ by maximal (under inclusion) such topological balls. Let $f^{n_i}$ be the corresponding iterates. 

Defining $\psi : \bigcup U_i \to U$ by $\psi_{|U_i} := f^{n_i}_{|U_i}$, we obtain an expanding induced Markov map for $(f,m)$. Let us check that it has integrable return time:

Firstly we remark that for each $i$, $n_i = \inf\{m(y) : y \in A \cap \Pi^{-1} U_i\}$ since the $U_i$ were maximal under inclusion.  
Thus, writing $n(x) := n_i$ if $x \in U_i$,
\begin{eqnarray*}
\sum_i n_i m(U_i) = \int_U n(x) dm 
&\leq& \delta^{-1} \int_U n(x) d\Pi_*\omu_A \\
&=&
 \delta^{-1} \int_A n(\Pi y) d\omu_A \\
&\leq& 
\delta^{-1} \int_A m(y) d\omu_A \\
&=&
\delta^{-1} \int_A m(y) d\omu \\
&=& \delta^{-1} N < \infty,
\end{eqnarray*}
as required, completing the proof (modulo the following variant of the folklore lemma) since $\mu$ is unique.
\eprf
\begin{lem}\label{lem:markovac}
Let $\psi : \bigcup_i U_i \to U$ be an expanding Markov map for $(f,m)$.
Then there is an absolutely continuous $\psi$-invariant probability measure $\nu$ on $U$. The density of $\nu$ with respect to $m$ is bounded away from zero and infinity.
\end{lem}
\beginpf By the Koebe principle, there is a uniform bound on the distortion of $\psi^n$ on connected components of its domain. We now want to show that $S_N \phi$ does not vary much either on a connected component $W$ of the domain of $\psi^n$, where $\psi^n = f^N$ on $W$. 

Recall that, by the definition of $\psi$, 
 there exists $C, \delta > 0$ such that, for each $i$ and all $j\leq n_i $, 
    $$
    |Df^j(x)| > C e^{\delta j}
    $$
    for all  $x \in  f^{n_i - j}(U_i)$ where $n_i$ is given by $\psi_{|U_i}= f^{n_i}.$ 

We wish to show that there are constants $C_0, \delta_0 >0$ such that $|f^j(W)| \leq C_0\exp((j-N)\delta_0)$ for all $j\leq N$. For this, using the previous paragraph, it suffices to show that there 
is a  constant $\delta_1$ such that $|W| \leq  \exp(-N\delta_1)$ for all such $W, N$. 
Let $M>0$ be small enough that $C^M  \exp(\delta) >1$. Then, considering the cases $n\geq NM$ and $n< NM$ we deduce, for all $x \in W$, 
$$
|D\psi^n(x)| = |Df^N(x)| > \min\left(2^{NM}, C^{NM} \exp(N\delta)\right)
 $$ 
 (since $|D\psi|>2$ everywhere). 
  We can set $\delta_1 := \log\min(2^M, C^M \exp(\delta))$. 

  Now, by H\"older continuity of $\phi$, there are $C', \epsilon>0$ such that for all $x,x' \in \cbar$, $|\phi(x) - \phi(x')| \leq C' |x-x'|^\epsilon$.  We deduce that if $x,x' \in W$ then 
  $$ 
  S_N \phi(x) - S_N\phi(x') 
  \leq \sum_{j=0}^\infty C' (C_0 \exp(-j\delta_0))^\epsilon < C_2$$
  for some constant $C_2 < \infty$, and this independently of $n$ and $W$. 
 
  This completes the first step. The rest of the proof is standard. 

 Since $m(U) = m(\psi^{-1}(U))$, for all $i\geq 1$, $\nu_i := \frac{1}{m(U)}(\psi^{-i})^* m$ is a probability measure. By the distortion bound and the bound on $
  S_N \phi(x) - S_N\phi(x') 
  $, there is a constant $C_3$ such that 
$$
C_3^{-1} \frac{m(A)}{m(B)} \leq \frac{\nu_i(A)}{\nu_i(B)} \leq C_3 \frac{m(A)}{m(B)}
$$
for all measurable sets $A,B \subset U$. The same relation holds with $\nu_i$ replaced by any weak limit $\nu$ of the sequence of the measures 
$$
\left( \frac{1}{n} \sum_{i=1}^n \nu_i\right)_{n=1}^\infty.
$$
Thus $\nu$ is a probability measure which is absolutely continuous with respect to $m$  with density bounded between $C_3^{-1}$ and $C_4$ (to see this, set $B = U$). Clearly $\nu$ is invariant. Ergodicity and uniqueness are easy.
\eprf 




\bibliography{references}
\bibliographystyle{plain}

\end{document}